\documentclass[12pt]{amsart}
\usepackage{amsfonts,amssymb,amsmath}
\usepackage[russian,english]{babel}
\usepackage[cp1251]{inputenc}
\usepackage{amsfonts,amssymb}
\usepackage{geometry}

\begin{document}
\selectlanguage{russian}

\title[Регуляризация двучленных дифференциальных уравнений]{Регуляризация двучленных дифференциальных уравнений с сингулярным коэффициентом}

\author{А.\,С.~Горюнов, В.\,А.~Михайлец}

\address{Институт математики НАН Украины, Киев}

\email[Андрей Горюнов]{goriunov@imath.kiev.ua} %Andrii Goriunov
\email[Владимир Михайлец]{mikhailets@imath.kiev.ua} %Vladimir Mikhailets


\begin{abstract}
 We propose a regularization of the formal differential expression of order $m \geqslant 3$
$$
l(y) = i^my^{(m)}(t) + q(t)y(t), \,t \in (a, b),
$$
applying quasi-derivatives.
The distribution coefficient $q$ is supposed to have an antiderivative $Q \in L\left([a,b];\mathbb{C}\right)$.
For the symmetric case ($Q = \overline{Q}$) self-adjoint and maximal dissipative extensions of the minimal operator and its
generalized resolvents are described.
The resolvent approximation with resrect to the norm of the considered operators is also investigated.

The case $m = 2$ for $Q \in L_2\left([a, b];\mathbb{C}\right)$
was investigated earlier.

%\vskip 1cm
%
%Предлагается регуляризация формального
%дифференциального выражения порядка $m \geqslant 3$
%$$
%l(y) = i^my^{(m)}(t) + q(t)y(t), \,t \in (a, b),
%$$
%посредством квазипроизводных. Предполагается, что коэффициент-распределение $q$ имеет первообразную $Q \in
%L\left([a,b];\mathbb{C}\right)$. В симметрическом случае ($Q =
%\overline{Q}$) описаны самосопряженные, максимальные диссипативные
%расширения минимального оператора и его обобщенные резольвенты.
%Изучена сходимость резольвент рассмотренных операторов по норме.
%
%Случай $m = 2$ при $Q \in L_2\left([a, b];\mathbb{C}\right)$
%исследован ранее.
\end{abstract}


\maketitle

%\noindent
%%%%%%%%%%%%%%%%%%%%%%%%%%%%%%%%%%%%%%%%%%%%%%%%%%%%%%%%%%%%%%%%%%%%%%%%%%%%%%%%%

%\textbf{1. Введение.}

\section{Введение.}
Рассмотрим на конечном интервале $\mathcal{J}~:=~(a, b)$ формальное дифференциальное выражение
порядка $m$
\begin{equation}\label{vyraz}
l(y) = i^m y^{(m)}(t) + q(t)y(t), \quad t \in \mathcal{J}.
\end{equation}
Если $m = 2$ и коэффициент $q \in L\left(\mathcal{J};\mathbb{R}\right)$, то дифференциальное уравнение $l(y) = f$ является
классическим уравнением Штурма--Лиувилля и изучено весьма полно.
Современное изложение этой теории приведено во многих монографиях (см. [\ref{Z}] и приведенные там ссылки). Как выяснилось после работы [\ref{S-Sh}], многие положения этой теории распространяются на существенно более общий случай
\begin{equation}\label{S-Sh cond}
q = Q', \quad Q \in L_2\left(\mathcal{J}; \mathbb{C}\right),
\end{equation}
где производная понимается в смысле обобщенных функций.
В частности, это относится к физически содержательному случаю, когда $q$ является мерой Радона на
$\overline{\mathcal{J}}$ либо имеет неинтегрируемые точечные особенности.
Подобные операторы задолго до этого возникали в различных задачах математической физики
и исследовались очень многими авторами,
главным образом, при помощи средств теории операторов (см. [\ref{Albeverio}] и ссылки там).

В связи с этим представляет интерес задача о регуляризации дифференциального выражения (\ref{vyraz}) с сингулярным
коэффициентом $q \notin L\left(\mathcal{J};\mathbb{C}\right)$ при произвольном значении $m > 2$.
Ее решению посредством специально подобранных квазипроизводных и посвящена данная работа.
При этом условие (\ref{S-Sh cond}) нам удалось ослабить до нижеследующего:
\begin{equation}\label{GM cond}
q = Q', \quad Q \in L\left(\mathcal{J}; \mathbb{C}\right) =: L_1.
\end{equation}

Случай общего выражения Штурма--Лиувилля
$$
l(y) = -(p(t)y')' + q(t)y, \quad t \in \mathcal{J}
$$
с сингулярными коэффициентами
$$
q = Q', \quad 1/p,\, Q/p,\, Q^2/p \in L_1
$$
с аналогичных позиций исследован авторами ранее в работе [\ref{GM St-L}].

Работа структурирована следующим образом.

В разделе 2 мы вводим регуляризацию формального дифференциального выражения (\ref{vyraz}) в предположении (\ref{GM cond})
и определяем соответствующие максимальный и минимальный операторы в гильбертовом пространстве
$L_2\left(\mathcal{J}; \mathbb{C}\right) =: L_2.$

В разделе 3 найдены достаточные условия равномерной резольвентной аппроксимации
расширений построенного минимального оператора $L_{\text{min}}$ семейством операторов того же класса.

В разделе 4 в предположении симметричности минимального оператора описываются все его самосопряженные,
максимальные диссипативные и максимальные аккумулятивные расширения в терминах однородных граничных условий
канонического вида.
Эти расширения параметризуются соответственно унитарными операторами и сжатиями в $\mathbb{C}^m$.
Такая параметризация является биективной и непрерывной.

В разделе 5 описываются все обобщенные резольвенты минимального оператора во всей комплексной плоскости.
%%%%%%%%%%%%%%%%%%%%%%%%%%%%%%%%%%%%%%%%%%
\section{Регуляризация сингулярного выражения.}
%\textbf{2. Регуляризация сингулярного выражения.}

 Рассмотрим формальное дифференциальное выражение (\ref{vyraz})
порядка $m \geq 3$ при условиях (\ref{GM cond}).

Введем последовательно квазипроизводные:
\begin{align*}
&D^{[k]}y(t) := y^{(k)}(t),\quad k=\overline{0,m-2},\\
&D^{[m-1]}y(t) := y^{(m-1)}(t)+i^{-m}Q(t)y(t),\\
&D^{[m]}y(t) := (D^{[m-1]}y(t))'-i^{-m}Q(t)D^{[1]}y(t).
\end{align*}

В условиях (\ref{GM cond}) они являются квазипроизводными по Шину--Зеттлу (см. [\ref{EM}], Sec. 1).

Поэтому формальное выражение (\ref{vyraz}) можно корректно определить как квазидифференциальное выражение Шина--Зеттла
$$l[y] := i^m D^{[m]}.$$

\textbf{Определение 1.}
\emph{Решение задачи Коши для резольвентного уравнения}
\begin{equation}\label{cauchy pr 1}
l[y] - \lambda y = f \in L_2, \qquad
 (D^{[k]}y)(c) = \alpha_k, \quad k=\overline{0,m-1},
\end{equation}
\emph{где} $c \in \overline{\mathcal{J}}$ \emph{и} $\alpha_k \in \mathbb{C}, \, k=\overline{0,m-1}$, \emph{определяется как первая компонента
решения задачи Коши для соответствующей системы дифференциальных уравнений первого порядка}
\begin{equation}\label{cauchy pr 2}
w'(t)=A_\lambda(t)w(t) + \varphi(t),
\qquad
w(c) = (\alpha_0, \alpha_1, ... \alpha_{m - 1})
\end{equation}
\emph{где вектор-функция} $w(t) := (D^{[0]}y(t), D^{[1]}y(t), ..., D^{[m - 1]}y(t)),$
\emph{квадратная матрица-функция}
\begin{equation}\label{A matrix}
A_\lambda(t):=\left (
\begin{array}{ccccccc}
0&1&0&\ldots &0&0 \\
0&0&1&\ldots &0&0 \\
\vdots &\vdots &\vdots &\ddots &\vdots &\vdots  \\
0&0&0&\ldots &1&0 \\
-i^{-m}Q(t)&0&0&\ldots &0&1\\
i^{-m}\lambda&\ i^{-m}Q(t)&0&\ldots &0 & 0
\end{array}\right) \in L_1^{m\times m},
\end{equation}
\emph{а вектор-функция} $\varphi(t) := (0, 0, ..., 0, i^{-m}f(t)) \in L_1^m$.

\textbf{Лемма 1.}
\emph{Задача Коши} (\ref{cauchy pr 1}) \emph{при условии} (\ref{GM cond}) \emph{имеет решение на $\overline{\mathcal{J}}$.
Оно единственно.}

\emph{\textbf{Доказательство.}}
Задача (\ref{cauchy pr 2}) при $A_\lambda(\cdot) \in L_1^{m\times m}$ имеет и притом единственное решение при
каждом $c \in \overline{\mathcal{J}}$
и $(\alpha_0, \alpha_2, ... \alpha_{m - 1}) \in \mathbb{C}^m$ в силу теоремы 1.2.1 монографии [\ref{Z}].
Поэтому утверждение леммы следует из определения 1 и указанной теоремы.

Введенное квазидифференциальное выражение $l[y]$  порождает в гильбертовом пространстве $L_2$
(см. [\ref{Zettl}, \ref{EM}])
\emph{максимальный} квазидифференциальный оператор
$$
L_{\text{max}}:y  \to l[y],
$$
%\begin{align*}
%\text{Dom}&(L_{\text{max}}) :=\\
%:= &\left\{y \in L_2 \left| D^{[k]}y \in AC(\overline{\mathcal{J}}, \mathbb{C}),\, k=\overline{0,m-1},\,
% D^{[m]} y \in L_2\right.\right\}.
%\end{align*}
$$
\text{\text{Dom}}(L_{\text{max}}) = \left\{y \left| D^{[k]}y \in AC(\overline{\mathcal{J}}, \mathbb{C}), k=\overline{0,m-1},\,
 D^{[m]} y \in L_2\right.\right\}.
 $$
\noindent\emph{Минимальный} квазидифференциальный оператор
определяется как сужение оператора $L_{\text{max}}$ на линейное
многообразие $\,\text{\text{Dom}}(L_{\text{min}}) :=$
%\begin{align*}
%\text{Dom}&(L_{\text{min}}) :=\\
%:= & \left\{y \in  \text{Dom}(L_{\text{max}}) \left| D^{[k]}y(a) = D^{[k]}y(b) = 0,\,
%k = \overline{0,m - 1}\right.\right\}.
%\end{align*}
$$
:= \left\{y \in  \text{\text{Dom}}(L_{\text{max}}) \left| D^{[k]}y(a) = D^{[k]}y(b) = 0,\,
k = \overline{0,m - 1}\right.\right\}.
$$

\textbf{Лемма 2.}
\emph{Если в равенствах} (\ref{vyraz}) \emph{и} (\ref{GM cond}) \emph{заменить выбранную первообразную $Q$ произвольной}
 $$\widetilde{Q} := Q + c, \quad c \in \mathbb{C},$$
\emph{то операторы} $L_{\text{max}},$ $L_{\text{min}}$ \emph{не изменятся.}

\textbf{\emph{Доказательство.}}
Покажем, что оператор $L_{\text{max}} =$ $= L_{\text{max}}(Q)$, совпадает с оператором $\widetilde{L}_{\text{max}} = L_{\text{max}}(\widetilde{Q})$.
Обозначим через $ \widetilde{D}^{[0]}y, \widetilde{D}^{[1]}y, ..., \widetilde{D}^{[m]}y$ квазипроизводные,
соответствующие отличной от $Q$ первообразной $\widetilde{Q}$.

Пусть $y \in \text{Dom}(L_{\text{max}})$.
Прямым подсчетом находим, что
\begin{align*}
&\widetilde{D}^{[0]}y = D^{[0]}y \in AC(\overline{\mathcal{J}}, \mathbb{C}), \\
&\ldots\ldots\ldots\ldots\ldots,\\
&\widetilde{D}^{[m - 2]}y = D^{[m - 2]}y \in AC(\overline{\mathcal{J}}, \mathbb{C}),\\
&\widetilde{D}^{[m - 1]}y = D^{[m - 1]}y + i^{-m}c\widetilde{D}^{[0]}y \in AC(\overline{\mathcal{J}}, \mathbb{C}),\\
&\widetilde{D}^{[m]}y = D^{[m]}y \in L_2.
\end{align*}
Это означает, что
\begin{align*}
\text{Dom}(L&_{\text{max}}) \subset \text{Dom}(\widetilde{L}_{\text{max}}) =\\
= &\left\{y \left| \widetilde{D}^{[k]}y \in AC(\overline{\mathcal{J}}, \mathbb{C}), k=\overline{0,m-1},
\widetilde{D}^{[m]} y \in L_2\right.\right\}.
\end{align*}

\noindent Аналогично показывается, что $\text{Dom}(L_{\text{max}}) \supset \text{Dom}(\widetilde{L}_{\text{max}})$.

\noindent Наконец,
$$\widetilde{L}_{\text{max}}y = i^m\widetilde{D}^{[m]}y = i^mD^{[m]}y = L_{\text{max}}y, \quad y \in \text{Dom}(L_{\text{max}}).$$

Покажем теперь, что $\widetilde{L}_{\text{min}} = L_{\text{min}}$.

\noindent Пусть $y \in \text{Dom}(L_{\text{min}})$.
Тогда
\begin{align*}
&\widetilde{D}^{[0]}y(a) = D^{[0]}y(a) = 0, \\
&\widetilde{D}^{[0]}y(b) = D^{[0]}y(b) = 0,\\
&\ldots\ldots\ldots\ldots\ldots\ldots\ldots,\\
&\widetilde{D}^{[m - 2]}y(a) = D^{[m - 2]}y(a) = 0, \\
&\widetilde{D}^{[m - 2]}y(b) = D^{[m - 2]}y(b) = 0,\\
&\widetilde{D}^{[m - 1]}y(a) = D^{[m - 1]}y(a) + i^{-m}c\widetilde{D}^{[0]}y(a) = 0 + 0 = 0, \\
&\widetilde{D}^{[m - 1]}y(b) = D^{[m - 1]}y(b) + i^{-m}c\widetilde{D}^{[0]}y(b) = 0 + 0 = 0.
\end{align*}
Это означает, что
$\text{Dom}(L_{\text{min}}) \subset \text{Dom}(\widetilde{L}_{\text{min}}).$
Аналогично устанавливается, что ${\text{Dom}(L_{\text{min}}) \supset \text{Dom}(\widetilde{L}_{\text{min}})}$.

Поскольку $\widetilde{L}_{\text{min}}y = \widetilde{L}_{\text{max}}y = L_{\text{max}}y = L_{\text{min}}y\,$ на функциях $y \in \text{Dom}(L_{\text{min}})$,
то лемма доказана.

Рассмотрим наряду с (\ref{vyraz}) формально сопряженное дифференциальное выражение
$$l^+(y) = i^m y^{(m)}(t) + \overline{q}(t)y(t),$$
где черта обозначает комплексное сопряжение.
Обозначим через $L^+_{\text{max}}$ и $L^+_{\text{min}}$ порождаемые им максимальный и минимальный операторы в пространстве $L_2$.
Тогда из результатов монографии [\ref{EM}] для общих квазидифференциальных выражений Шина--Зеттла
и приведенного нами выше следует

\textbf{Теорема 1.}
\emph{Операторы} $L_{\text{min}}$, $L^+_{\text{min}}$, $L_{\text{max}}$, $L^+_{\text{max}}$ \emph{являются плотно заданными и замкнутыми в пространстве $L_2$,}
$$L_{\text{min}}^* = L^+_{\text{max}},\quad L_{\text{max}}^* = L^+_{\text{min}}.$$
\emph{Если функция $q$ вещественнозначна, то оператор} $L_{\text{min}} = L^+_{\text{min}}$ \emph{является симметрическим
с индексом дефекта} $\left({m,m} \right)$ \emph{и}
$$L_{\text{min}}^* = L_{\text{max}},\quad L_{\text{max}}^* = L_{\text{min}}.$$

%%%%%%%%%%%%%%%%%%%%%%%%%%%%%%%%%%%%%%%%%%%%%%%%%%%%%%%%

\section{Аппроксимация резольвенты.}
%\textbf{3. Аппроксимация резольвенты.}

Рассмотрим семейство квазидифференциальных выражений вида (\ref{vyraz}) $l_\varepsilon[y]$ с коэффициентами
$q_\varepsilon = Q'_\varepsilon \in L_1$, $\varepsilon \in [0, \varepsilon_0]$.
Соответствующие им квазипроизводные будем обозначать через
$D_\varepsilon^{[0]}y, D_\varepsilon^{[1]}y, ..., D_\varepsilon^{[m]}y$.

В гильбертовом пространстве $L_2$ с нормой $\|\cdot\|_2$ такие выражения при каждом $\varepsilon$
порождают операторы $L^\varepsilon_{\text{min}} $, $L^\varepsilon_{\text{max}}$.
Пусть матрицы $\alpha(\varepsilon),\beta(\varepsilon)\in \mathbb{C}^{m\times m},$
а векторы
%\begin{equation}\label{Y_a}
$$
\mathcal{Y}_\varepsilon(a):=\{D^{[0]}_\varepsilon y(a),D^{[1]}_\varepsilon y(a),\ldots,
D^{[m-1]}_\varepsilon y(a)\}\in\mathbb{C}^m,
$$
%\end{equation}
%\begin{equation}\label{Y_b}
$$
\mathcal{Y}_\varepsilon(b):=\{D^{[0]}_\varepsilon y(b),D^{[1]}_\varepsilon y(b),\ldots,
D^{[m-1]}_\varepsilon y(b)\}\in\mathbb{C}^m.
$$

Зададим для каждого фиксированного значения $\varepsilon$ операторы
$$L_\varepsilon y = l_\varepsilon[y],$$
$$
\text{Dom}(L_\varepsilon) = \left\{\left.y \in \text{Dom}\left(L^\varepsilon_{\text{max}}\right)\right|
\alpha(\varepsilon)\mathcal{Y}_\varepsilon(a)+\beta(\varepsilon)\mathcal{Y}_\varepsilon(b)=0\right\}.
$$

Очевидно, что
$$L^\varepsilon_{\text{min}} \subset L_\varepsilon \subset L^\varepsilon_{\text{max}}, \quad \varepsilon \in [0, \varepsilon_0].$$

Будем обозначать через $\rho(L)$ резольвентное множество оператора $L$.
Напомним, что операторы $L_\varepsilon$ сходятся при ${\varepsilon\rightarrow 0+}$ к оператору $L_0$
в смысле равномерной резольвентной сходимости,
$L_\varepsilon \stackrel{R}{\rightarrow} L_0$,
если существует $\mu \in \mathbb{C}$
такое, что $\mu \in \rho(L_0)$, $\mu \in \rho(L_\varepsilon)$ для достаточно малых $\varepsilon$ и
$$\|(L_\varepsilon - \mu)^{-1} - (L_0 - \mu)^{-1}\| \rightarrow 0, \quad \varepsilon \rightarrow 0+.$$
Это определение не зависит от выбора $\mu \in \rho(L_0)$ [\ref{K}].

Введем обозначение $c^\vee (t):= \int\limits_a ^t c(s)ds.$

Основным результатом этого раздела является

\textbf{Теорема 2.}
\emph{Пусть $\rho(L_0)$ непусто и при $\varepsilon\rightarrow 0+$ выполняются условия:}

 $1) \quad
\|(Q_\varepsilon - Q_0)^\vee\|_C \rightarrow 0;$

$2) \quad \alpha(\varepsilon){\longrightarrow}\alpha(0),\quad
\beta(\varepsilon){\longrightarrow}\beta(0).$

\emph{Тогда $L_\varepsilon \stackrel{R}{\rightarrow} L_0$.}

\textbf{\emph{Замечание 1.}}
Условие $\|Q_\varepsilon - Q_0\|_1 \rightarrow 0, \, \varepsilon\rightarrow 0+$, очевидно, достаточно для выполнения условия
$1)$.

Доказательство теоремы 2 основывается на одном вспомогательном результате.

Следуя работам [\ref{MR1}, \ref{MR2}], введем

\textbf{Определение 2.}
\emph{Обозначим через $\mathcal{M}^{n}(\mathcal{J})=:\mathcal{M}^{n},$ $n \in \mathbb{N}$
класс всех параметризованных числом $\varepsilon $ матриц-функций}
$$R(\cdot;\varepsilon):[0,\varepsilon_0]\rightarrow L_1 ^{n\times n}, $$
\emph{для которых решение задачи Коши}
$$ Z'(t;\varepsilon)= R( t;\varepsilon)Z(t;\varepsilon), \quad Z(a;\varepsilon) = I_n$$
\emph{удовлетворяет предельному соотношению}
$$\lim\limits_{\varepsilon \rightarrow 0+} \|Z(\cdot;\varepsilon) - I_n\|_C =0,$$
\emph{где $\|\cdot\|_C$ -- $\sup$-норма.}

В работе [\ref{MR2}] установлена следующая общая теорема.

\textbf{Теорема 3.}
\emph{Пусть для краевой задачи}
\begin{equation}\label{bound probl 1}
  y'(t;\varepsilon)=A(t;\varepsilon)y(t;\varepsilon)+f(t;\varepsilon),\quad
     t \in\mathcal{J}, \quad \varepsilon \in [0, \varepsilon_0]
\end{equation}
\begin{equation}\label{bound probl 2}
    U_{\varepsilon}y(\cdot;\varepsilon)= 0,
\end{equation}
\emph{где матрицы-функции $A(\cdot,\varepsilon) \in  L_{1}^{n\times n}$,
вектор-функции ${f(\cdot,\varepsilon) \in L_{1}^{n}}$,\,
а линейные непрерывные операторы}
$$U_{\varepsilon}:C(\overline{\mathcal{J}};\mathbb{C}^{n})  \rightarrow\mathbb{C}^{n}, \, n \in \mathbb{N},$$
\emph{выполнены условия:}
\begin{align*}
 1) \quad &\text{\emph{однородная предельная краевая задача} ~}  (\ref{bound probl 1}), (\ref{bound probl 2})
 \text{~с~} \varepsilon = 0  \text{~\emph{и}~}  \\
  &f(\cdot;0) \equiv 0 \text{~\emph{имеет только тривиальное решение;}}\\
  2)\quad &A(\cdot;\varepsilon)-A(\cdot;0)\in \mathcal{M}^n;\\
  3)\quad &\|U_{\varepsilon} - U_{0}\|\rightarrow 0,\quad \varepsilon\rightarrow 0+.
\end{align*}
\emph{Тогда для достаточно малых $\varepsilon$
 существуют матрицы Грина $G(t, s; \varepsilon)$ задач} (\ref{bound probl 1}), (\ref{bound probl 2})
\emph{и на квадрате $\mathcal{J}\times \mathcal{J}$}
\begin{equation}\label{G}
 \|G(\cdot,\cdot;\varepsilon)-G(\cdot,\cdot;0)\|_\infty \rightarrow 0,\quad
\varepsilon\rightarrow 0+,
\end{equation}
\emph{где $\|\cdot\|_\infty$ -- норма в пространстве $L_\infty$.}

\textbf{\emph{Замечание 2.}}
Условие 3) теоремы 3 нельзя заменить более слабым условием сильной сходимости операторов
$U_{\varepsilon}\stackrel{s}\rightarrow U_{0}$ [\ref{MR2}].
Однако, как нетрудно убедиться, для двухточечных краевых операторов
$$U_\varepsilon y := B_1(\varepsilon) y(a) + B_2(\varepsilon) y(b),
\quad B_k(\varepsilon) \in \mathbb{C}^{n\times n},\quad k \in \{1,2\},$$
как условие сильной, так и условие равномерной сходимости равносильны тому, что
$${\|B_k(\varepsilon) - B_k(0)\| \rightarrow 0,}\quad \varepsilon\rightarrow 0+, \quad k \in \{1,2\}.$$

Приведенное определение класса $\mathcal{M}^n$ не является конструктивным.
Имеются различные достаточные условия принадлежности матричной функции $R(\cdot;\varepsilon)$ классу
$\mathcal{M}^{n}$.
В частности, из результатов работы А. Ю. Левина [\ref{Levin}] следует

\textbf{Лемма 3.}
\emph{Пусть $R(\cdot;\varepsilon):[0,\varepsilon_0]\rightarrow L_1 ^{n\times n}.$
Если при  $\varepsilon\rightarrow 0+$
выполнено одно из четырех (неэквивалентных между собой) условий:}
\begin{align*}
(\alpha)& \quad \|R(\cdot;\varepsilon)\|_1 = O(1),\quad\quad\quad\quad\quad\quad\quad\quad\quad\quad
\quad\quad\quad\quad\quad\quad\quad\quad\\
(\beta)& \quad\|R^\vee(\cdot;\varepsilon)R(\cdot;\varepsilon)\|_1 \rightarrow 0,\\
(\gamma)& \quad\|R(\cdot;\varepsilon)R^\vee(\cdot;\varepsilon)\|_1\rightarrow 0,\\
(\Delta)& \quad\|R^\vee(\cdot;\varepsilon)R(\cdot;\varepsilon) -
           R(\cdot;\varepsilon)R^\vee(\cdot;\varepsilon)\|_1 \rightarrow 0,
\end{align*}
\emph{то условие
$\|R^\vee (\cdot;\varepsilon)\|_C\rightarrow 0,\, \varepsilon\rightarrow 0+$
равносильно включению
$R(\cdot;\varepsilon) \in \mathcal{M}^{n}$.}

Следующее утверждение позволит нам редуцировать теорему 2 к теореме 3.

\textbf{Лемма 4.}
\emph{Функция $y(t)$ является решением краевой задачи}
\begin{equation}\label{D^m}
l_\varepsilon[y](t)= f(t;\varepsilon) \in L_2 ,\quad\varepsilon\in [0,\varepsilon_0],
\end{equation}
\begin{equation}\label{alpha+beta}
  \alpha(\varepsilon)\mathcal{Y}_\varepsilon(a)+
  \beta(\varepsilon)\mathcal{Y}_\varepsilon(b)=0
\end{equation}
\emph{тогда и только тогда, когда вектор-функция}
$$w(t) = (D_\varepsilon^{[0]}y(t), D_\varepsilon^{[1]}y(t), ..., D_\varepsilon^{[m - 1]}y(t))$$
\emph{является решением краевой задачи}
\begin{equation}\label{diff eq}
w'(t)=A(t;\varepsilon)w(t) + \varphi(t;\varepsilon),
\end{equation}
\begin{equation}\label{diff alpha+beta}
\alpha(\varepsilon)w(a)+
  \beta(\varepsilon)w(b)=0,
  \end{equation}
\emph{где квадратная матрица-функция}
\begin{equation}\label{A matrix}
A(\cdot;\varepsilon):=\left (
\begin{array}{ccccccc}
0&1&0&\ldots &0&0 \\
0&0&1&\ldots &0&0 \\
\vdots &\vdots &\vdots &\ddots &\vdots &\vdots  \\
0&0&0&\ldots &1&0 \\
-i^{-m} Q(\cdot;\varepsilon)&0&0&\ldots &0&1\\
0&\ i^{-m}Q(\cdot;\varepsilon)&0&\ldots &0 & 0
\end{array}\right) \in L_1^{m\times m},
\end{equation}
\emph{а} $\varphi(\cdot;\varepsilon) := (0, 0, ..., 0, i^{-m}f(\cdot;\varepsilon)) \in L_1^m$.

\textbf{\emph{Доказательство.}}
 Рассмотрим систему уравнений
$$\left\{
\begin{array}{l}
    ( D^{[0]}_\varepsilon y(t))' = D^{[1]}_\varepsilon y(t) \\
    ( D^{[1]}_\varepsilon y(t))' = D^{[2]}_\varepsilon y(t) \\
    \ldots\ldots\ldots\ldots\ldots\ldots\ldots\ldots\ldots, \\
    ( D^{[m - 3]}_\varepsilon y(t))' = D^{[m - 2]}_\varepsilon y(t) \\
    ( D^{[m - 2]}_\varepsilon y(t))' = -i^{-m} Q_\varepsilon(t)D^{[0]}_\varepsilon y(t) +
          D^{[m - 1]}_\varepsilon y(t) \\
    ( D^{[m - 1]}_\varepsilon y(t))' = i^{-m}Q_\varepsilon(t)D^{[1]}_\varepsilon y(t) +
           i^{-m} f(t; \varepsilon)\\
\end{array}
\right.$$

Если $y(\cdot)$ -- решение уравнения (\ref{D^m}), то из определения
квазипроизводных следует, что $y(\cdot)$ есть решение этой системы.
С другой стороны, положив
$w(t) = (D^{[0]}_\varepsilon y(t), D^{[1]}_\varepsilon y(t), ...,  D^{[m - 1]}_\varepsilon y(t))$
и $\varphi(t;\varepsilon)~=$ $=~(0,0, ..., 0,  i^{-m}f(t;\varepsilon))$, данную систему можно записать в виде уравнения
(\ref{diff eq}).

Учитывая, что $\mathcal{Y}_\varepsilon(a) = w(a)$,
$\mathcal{Y}_\varepsilon(b) = w(b)$, легко видеть, что краевые
условия (\ref{alpha+beta}) эквивалентны краевым условиям (\ref{diff
alpha+beta}).

В силу леммы 4 из предположения
\begin{itemize}
\item [$(\mathcal{E})$]
\emph{Однородная краевая задача} $$D^{[m]}_0 y(t)=0, \quad \alpha(0)\mathcal{Y}_0(a)+ \beta(0)\mathcal{Y}_0(b) = 0$$
\emph{имеет только тривиальное решение}
\end{itemize}
следует, что однородная краевая задача
$$w'(t)=A(t;\varepsilon)w(t), \quad \alpha(\varepsilon)w(a) + \beta(\varepsilon)w(b)=0$$
также имеет только тривиальное решение.

\textbf{Лемма 5.}
\emph{Пусть для задачи} (\ref{diff eq}), (\ref{diff alpha+beta}) \emph{при достаточно малых $\varepsilon$ существует матрица Грина}
$$G(t,s,\varepsilon)=(g_{ij}(t,s))_{i,j=1}^m \in L_\infty^{m\times m}$$
\emph{Тогда существует функция Грина $\Gamma(t,s;\varepsilon)$ полуоднородной краевой задачи} (\ref{D^m}), (\ref{alpha+beta})
\emph{и}
$$ \Gamma(t,s;\varepsilon) =  i^{-m} g_{1m}(t,s;\varepsilon) \qquad\mbox{\emph{п. в.}}$$

\textbf{\emph{Доказательство.}}
Согласно определению матрицы Грина, единственное решение задачи (\ref{diff eq}), (\ref{diff alpha+beta})
записывается в виде
$$w_\varepsilon(t)=\int\limits_a ^b G(t,s;\varepsilon)\varphi(s;\varepsilon) ds, \quad t\in \overline{\mathcal{J}}.$$

В силу леммы 4 последнее равенство можно переписать в виде
$$\left\{
\begin{array}{l}
    D^{[0]}_\varepsilon y_\varepsilon(t) = \int\limits_a^b g_{1m}(t,s;\varepsilon) i^{-m}f(s;\varepsilon)ds \\
    D^{[1]}_\varepsilon y_\varepsilon(t) = \int\limits_a^b g_{2m}(t,s;\varepsilon) i^{-m}f(s;\varepsilon)ds, \\
    \ldots\ldots\ldots\ldots\ldots\ldots\ldots\ldots\ldots\ldots\ldots\ldots, \\
    D^{[m]}_\varepsilon y_\varepsilon(t) = \int\limits_a^b g_{mm}(t,s;\varepsilon) i^{-m}f(s;\varepsilon)ds, \\
\end{array}
\right.$$ где $y_\varepsilon(\cdot)$ -- единственное решение задачи
(\ref{D^m}), (\ref{alpha+beta}). Отсюда следует утверждение леммы 5.

\textbf{Доказательство теоремы 2.}
В силу условия 1) теоремы 2 можно, не уменьшая общности, считать, что $0 \in \rho(L_0)$.

Покажем, что
$\sup\limits_{\|f\|_2 =1} \|L_\varepsilon^{-1}f - L_0^{-1}f\| \rightarrow 0$, $\varepsilon \rightarrow 0+$.

Уравнение $L_\varepsilon^{-1}f = y_\varepsilon$ эквивалентно тому,
что $L_\varepsilon y_\varepsilon = f$, т. е. $y_\varepsilon$
является решением задачи (\ref{D^m}), (\ref{alpha+beta}). При этом
из включения $0 \in \rho(L_0)$ следует, что выполняется
предположение ($\mathcal{E}$).

Обозначим $r(\cdot;\varepsilon) := i^{-m}Q(\cdot;\varepsilon) - i^{-m}Q(\cdot;0)$.
Тогда
$$A(\cdot;\varepsilon) - A(\cdot;0) = \left (
\begin{array}{ccccccc}
0&0&0&\ldots &0&0 \\
0&0&0&\ldots &0&0 \\
\vdots &\vdots &\vdots &\ddots &\vdots &\vdots  \\
-r(\cdot;\varepsilon)&0&0&\ldots &0&0\\
0&\ r(\cdot;\varepsilon)&0&\ldots &0 & 0
\end{array}\right),$$
$$\left(A(\cdot;\varepsilon) - A(\cdot;0)\right)^\vee = \left (
\begin{array}{ccccccc}
0&0&0&\ldots &0&0 \\
0&0&0&\ldots &0&0 \\
\vdots &\vdots &\vdots &\ddots &\vdots &\vdots  \\
-r^\vee(\cdot;\varepsilon)&0&0&\ldots &0&0\\
0&\ r^\vee(\cdot;\varepsilon)&0&\ldots &0 & 0
\end{array}\right),$$
где матрица-функция $A(\cdot;\varepsilon)$ задана формулой (\ref{A matrix}).

Легко видеть, что
\begin{align*}
\left(A(\cdot;\varepsilon) - A(\cdot;0)\right)&\left(A(\cdot;\varepsilon) - A(\cdot;0)\right)^\vee =\\
&= \left(A(\cdot;\varepsilon) - A(\cdot;0)\right)^\vee\left(A(\cdot;\varepsilon) - A(\cdot;0)\right).
\end{align*}
Поэтому матричная функция $A(\cdot;\varepsilon) - A(\cdot;0)$ при $m \geq 3$ удовлетворяет условию
$(\Delta)$ леммы 3.

Очевидно, что условие $\|\left(A(\cdot;\varepsilon)-A(\cdot;0)\right)^\vee\|_C\rightarrow 0,\, \varepsilon\rightarrow 0+$
эквивалентно условию $1)$ теоремы 2.
Поэтому из леммы 3 вытекает, что выполнены условия теоремы 3
для задачи (\ref{diff eq}), (\ref{diff alpha+beta}).

Это значит, что существуют матрицы Грина $G(t,s;\varepsilon)$ задач (\ref{diff eq}), (\ref{diff alpha+beta}),
и выполняется предельное соотношение (\ref{G}).
Учитывая лемму 5, это влечет предельное равенство
$$ \|\Gamma(\cdot,\cdot;\varepsilon)-\Gamma(\cdot,\cdot;0)\|_\infty \rightarrow 0,
\quad \varepsilon\rightarrow 0+.$$

Тогда

\parshape 3 0cm 15 cm 1.5 cm 10 cm 3 cm 10 cm \noindent$\|L_\varepsilon^{-1}- L_0^{-1}\| =$
$\sup\limits_{\|f\|_2=1} \|\int_a^b[\Gamma(t,s;\varepsilon) - \Gamma(t,s;0)]f(s)ds\|_2 \leq $\\
$\leq {(b - a)^{1/2}}{\sup\limits_{\|f\|_2=1}\|\int_a^b
|\Gamma(t,s;\varepsilon) - \Gamma(t,s;0)| |f(s)| ds\|_C} \leq$\\
$\leq (b - a)\|\Gamma(\cdot,\cdot;\varepsilon) - \Gamma(\cdot,\cdot;0)\|_\infty \rightarrow 0,
\quad \varepsilon \rightarrow 0+$,

\noindent что влечет утверждение теоремы 2.

%%%%%%%%%%%%%%%%%%%%%%%%%%%%%%%%%%%%%%%%%%%%%%%%%%%%%%%%

\section{ Расширения симметрического минимального оператора.}

%\textbf{4. Расширения симметрического минимального оператора.}

Всюду далее будем предполагать, что функции $q$ и, следовательно, $Q$ вещественнозначны.
Это условие влечет формальную самосопряженность выражения $l[y]$ (см. [\ref{EM}]) и, согласно теореме 1,
симметричность оператора $L_{\text{min}}$.
Поэтому содержателен вопрос об описании (при помощи однородных краевых условий)
некоторых классов (самосопряженных, максимальных диссипативных, максимальных аккумулятивных) расширений в
гильбертовом пространстве $L_2$ симметрического оператора $L_{\operatorname{min}}$.
Для ответа на него мы будем использовать понятие пространства граничных значений.

\textbf{Определение 3.}
\emph{Пусть $L$ -- замкнутый симметрический оператор в гильбертовом пространстве} $\mathcal{H}$ \emph{с равными
(конечными или бесконечными) дефектными числами.}
\emph{Тройка }$\left( {H,\Gamma _1 ,\Gamma _2 }\right)$\emph{, где} $H$\emph{ -- вспомогательное гильбертово пространство,
а} $\Gamma_1$, $\Gamma_2$ \emph{-- линейные отображения} $\text{Dom}(L^*)$\emph{ в }$H,$
\emph{называется пространством граничных значений (сокращенно ПГЗ) симметрического оператора $L$, если:}
%\begin{enumerate}
%\item
%для любых $f,g \in \text{Dom} \left( {L^*} \right)$
%  $\left( {L^ *  f,g} \right)_\mathcal{H} - \left( {f,L^ *  g}\right)_\mathcal{H} =
%  \left( {\Gamma _1 f,\Gamma _2 g} \right)_H  - \left( {\Gamma _2 f,\Gamma _1 g} \right)_H,$
%\item для любых векторов  $f_1, f_2 \in H$
%существует вектор $f\in \text{Dom} \left( {L^*} \right)$ такой, что
%$\Gamma _1 f = f_1,\,  \Gamma _2 f = f_2.$
%\end{enumerate}
%\begin{align*}
% 1) \quad &\text{для любых} f,g \in \text{Dom} \left( {L^*} \right)\\
%          &\left( {L^ *  f,g} \right)_\mathcal{H} - \left( {f,L^ *  g}\right)_\mathcal{H} =
%           \left( {\Gamma _1 f,\Gamma _2 g} \right)_H  - \left( {\Gamma _2 f,\Gamma _1 g} \right)_H,\\
%  2) \quad &\text{для любых векторов} f_1, f_2 \in H\\
%     &\text{существует вектор} f\in \text{Dom} \left( {L^*} \right) \text{такой, что}\\
%     &\Gamma _1 f = f_1,\,  \Gamma _2 f = f_2.
%\end{align*}

$1) \quad $\emph{для любых }$f,g \in \text{Dom} \left( {L^*} \right)$

  $\left( {L^ *  f,g} \right)_\mathcal{H} - \left( {f,L^ *  g}\right)_\mathcal{H} =
  \left( {\Gamma _1 f,\Gamma _2 g} \right)_H  - \left( {\Gamma _2 f,\Gamma _1 g} \right)_H,$

$2) \quad$\emph{ для любых векторов  }$f_1, f_2 \in H$
\emph{существует вектор} ${f\in \text{Dom} \left( {L^*} \right)}$ \emph{такой, что}
$\Gamma _1 f = f_1,\,  \Gamma _2 f = f_2.$

Из определения ПГЗ следует, что $ f \in \text{Dom} \left( {L} \right)$ тогда и только тогда,
когда $\Gamma_1f = \Gamma_2f = 0$.
ПГЗ существует для любого симметрического оператора с равными ненулевыми дефектными числами
(см. [\ref{Gorbachuk}] и приведенные там ссылки).
Оно всегда не единственно.

Следующий результат является ключевым для дальнейшего.

\textbf{Основная лемма.}
\emph{Пусть $\Gamma_1, \Gamma_2$ -- линейные отображения из} $\text{Dom}(L_{\text{max}})$ \emph{в} $\mathbb{C}^{m}$ \emph{такие, что:}

\noindent\emph{при} $m = 2n, n \geq 2,$
%\begin{equation} \label{PGZ 2n}
$$\Gamma_1y := i^{2n}
 \left(
   \begin{array}{c}
   - D^{[2n - 1]}y(a),\\
   ...,\\
   (-1)^nD^{[n]}y(a),\\
   D^{[2n - 1]}y(b),\\
   ...,\\
   (-1)^{n - 1}D^{[n]}y(b)
   \end{array}
 \right),
 \Gamma_2y :=
  \left(
   \begin{array}{c}
   D^{[0]}y(a),\\
   ...,\\
   D^{[n - 1]}y(a),\\
   D^{[0]}y(b),\\
   ...,\\
   D^{[n - 1]}y(b)
   \end{array}
  \right)$$
%\end{equation}
\emph{а при }$m = 2n + 1, n \in \mathbb{N},$

%\begin{equation} \label{PGZ 2n+1}
 %\begin{align}
 $$
  \Gamma_1y := i^{2n + 1}
   \left(
    \begin{array}{c}
    - D^{[2n]}y(a),\\
    ...,\\
    (-1)^{n} D^{[n + 1]}y(a),\\
    D^{[2n]}y(b),\\
    ....,\\
    (-1)^{n - 1}D^{[n + 1]}y(b),\\
    \alpha D^{[n]}y(b) + \beta D^{[n]}y(a)
    \end{array}
   \right),
   $$

   $$
  \quad\quad\Gamma_2y :=
   \left(
   \begin{array}{c}
   D^{[0]}y(a),\\
   ...,\\
   D^{[n -1]}y(a),\\
   D^{[0]}y(b),\\
   ...,\\
   D^{[n - 1]}y(b),\\
   \gamma D^{[n]}y(b) + \delta D^{[n]}y(a)
   \end{array}
  \right),
  $$
 %\end{align}
%\end{equation}
\noindent \emph{где числа} $\alpha = 1, \beta = 1, \gamma = \frac{(-1)^n}{2} + i, \delta = \frac{(-1)^{n + 1}}{2} + i.$

\emph{Тогда тройка} $(\mathbb{C}^{m}, \Gamma_1, \Gamma_2)$\emph{ является пространством граничных значений оператора }$L_{\text{min}}.$

\textbf{\emph{Замечание 3.}}
Приведенные значения коэффициентов можно заменить
произвольными наборами чисел, которые удовлетворяют системе:
\begin{equation}\label{PGZ coef}
\left\{
\begin{array}{l}
    \alpha\overline{\gamma} + \overline{\alpha}\gamma = (-1)^{n},\\
    \beta\overline{\delta} + \overline{\beta}\delta = (-1)^{n + 1}, \\
    \alpha\overline{\delta} + \overline{\beta}\gamma = 0, \\
    \beta\overline{\gamma}  + \overline{\alpha}\delta = 0, \\
    \alpha\delta - \beta\gamma \neq 0.\\
\end{array}
\right.
\end{equation}

Обозначим через $L_K$ сужение оператора $L_{\text{max}}$ на множество функций
${y(t) \in \text{Dom}(L_{\text{max}})}$, удовлетворяющих однородному краевому условию канонического вида
\begin{equation} \label{rozsh}
 \left( {K - I} \right)\Gamma _1 y + i\left( {K + I} \right)\Gamma _2 y = 0,
\end{equation}
где $K$ -- ограниченный оператор в гильбертовом пространстве $\mathbb{C}^{m}$.

Аналогично обозначим через $L^K$ сужение оператора $L_{\text{max}}$ на множество функций
${y(t) \in \text{Dom}(L_{\text{max}})}$, удовлетворяющих однородному краевому условию канонического вида
\begin{equation} \label{rozsh}
 \left( {K - I} \right)\Gamma _1 y - i\left( {K + I} \right)\Gamma _2 y = 0,
\end{equation}
где $K$ -- ограниченный оператор в гильбертовом пространстве $\mathbb{C}^{m}$.

В сочетании с результатами [\ref{Gorbachuk}]
Основная лемма влечет следующее описание самосопряженных расширений $L_{\text{min}}$:

\textbf{Теорема 4.}
\emph{Каждое $L_K$, где $K$ -- унитарный оператор в пространстве} $\mathbb{C}^{m}$,
\emph{является самосопряженным расширением оператора} $L_{\text{min}}.$
\emph{Обратно, для каждого самосопряженного расширения}
$\widetilde{L}$ \emph{оператора} $L_{\text{min}}$ \emph{найдется унитарный оператор $K$ такой,
что }$\widetilde{L} = L_K$.
\emph{Соответствие между унитарными операторами $\{K\}$ и расширениями }$\{\widetilde{L}\}$\emph{ биективно.}

\textbf{\emph{Замечание 4.}}
Из теоремы 2 и теоремы 4 вытекает, что отображение $K \to L_K$ является не только
биективным, но и непрерывным.
Более точно, если унитарные операторы $K_n$ сходятся по норме к оператору $K$, то
$$
\left\|\left(L_K - \lambda\right)^{-1} - \left(L_{K_n} - \lambda\right)^{-1}\right\| \rightarrow 0,
\,\, n \rightarrow \infty, \,\, \text{Im} \lambda \neq 0.
$$
При этом, поскольку множество унитарных операторов в конечномерном пространстве  $\mathbb{C}^m$ компактно
в метрике нормы оператора,
то верно и обратное утверждение, то есть отображение
$$K \to \left(L_K - \lambda\right)^{-1}, \, \text{Im} \lambda \neq 0$$
является при каждом фиксированном $\lambda \in \mathbb{C}\setminus\mathbb{R}$ гомеоморфизмом.

Напомним известное определение.

\textbf{Определение 4.}
\emph{Плотно заданный линейный оператор }$L$\emph{ в комплексном гильбертовом пространстве }$\mathcal{H}$ \emph{называют диссипативным (аккумулятивным),
если}
 $$\text{Im} \left( Lf, f \right)_\mathcal{H} \geq 0 \quad (\leq 0), \quad f \in \text{Dom} (L) $$
\emph{и максимальным диссипативным (максимальным аккумулятивным), если, кроме того, у оператора $L$ нет
нетривиальных диссипативных (аккумулятивных) расширений в пространстве} $\mathcal{H}.$

В частности, каждый симметрический оператор -- диссипативный и аккумулятивный, а
самосопряженный -- максимальный диссипативный и максимальный аккумулятивный одновременно.
Поэтому для симметрического квазидифференциального оператора $L_{\text{min}}$ можно поставить вопрос об описании
всех его максимальных диссипативных и максимальных аккумулятивных расширений.
Согласно теореме Р. Филлипса [\ref{Gorbachuk}, \ref{Phil}]
\emph{каждое диссипативное и каждое аккумулятивное расширение симметрического оператора является сужением его
сопряженного.}
Поэтому любое максимальное диссипативное или максимальное аккумулятивное расширение оператора $L_{\text{min}}$ является сужением оператора $L_{\text{max}}$.

Параметрическое описание всех максимальных диссипативных расширений
симметрического квазидифференциального оператора $L_{\text{min}}$ дает

\textbf{Теорема 5.}
\emph{Каждое} $L_K$\emph{, где} $K$ \emph{-- сжатие в пространстве }$\mathbb{C}^{m}$,
\emph{является максимально диссипативным расширением }$L_K$ \emph{оператора} $L_{\text{min}}.$
\emph{Обратно, для каждого максимального диссипативного расширения }$\widetilde{L}$ \emph{оператора} $L_{\text{min}}$
\emph{найдется сжатие }$K$\emph{ такое, что }$\widetilde{L} = L_K$.
\emph{Соответствие между сжатиями }$\{K\}$ \emph{и расширениями} $\{\widetilde{L}\}$ \emph{биективно.}

Параметрическое описание всех максимальных аккумулятивных расширений
симметрического квазидифференциального оператора $L_{\text{min}}$ дает

\textbf{Теорема 6.}
\emph{Каждое} $L^K$\emph{, где} $K$ \emph{-- сжатие в пространстве }$\mathbb{C}^{m}$,
\emph{является максимально аккумулятивным расширением }$L^K$ \emph{оператора} $L_{\text{min}}.$
\emph{Обратно, для каждого максимального аккумулятивного расширения }$\widetilde{L}$ \emph{оператора} $L_{\text{min}}$
\emph{найдется сжатие }$K$\emph{ такое, что }$\widetilde{L} = L^K$.
\emph{Соответствие между сжатиями }$\{K\}$ \emph{и расширениями} $\{\widetilde{L}\}$ \emph{биективно.}

%\emph{Для произвольного сжатия $K$ в пространстве $\mathbb{C}^{m}$ сужение $L_K$ оператора $L_{\operatorname{max}}$
%на множество функций, удовлетворяющих однородному краевому условию канонического вида}
%\begin{equation} \label{akk rozsh high}
% \left( {K - I} \right)\Gamma _{1} y - i\left( {K + I} \right)\Gamma _{2} y = 0,
%\end{equation}
%\emph{является максимально аккумулятивным расширением оператора $L_{\operatorname{min}}$.}
%\emph{Обратно, для каждого максимального аккумулятивного расширения $\widetilde{L}$ оператора $L_{\operatorname{min}}$
%найдется сжатие $K$ в пространстве $\mathbb{C}^{m}$ такое, что $\widetilde{L} = L_K$.}
%\emph{Соответствие между сжатиями }$\{K\}$ \emph{и расширениями} $\{\widetilde{L}\}$ \emph{биективно.}

\textbf{\emph{Замечание 5.}}
Отображения
%\\
%$\quad  \hfill K \rightarrow \left(L_K - \lambda\right)^{-1}, \, Im \lambda < 0$\hfill\quad\\
$$K \to \left(L_K - \lambda\right)^{-1}, \, Im \lambda < 0$$
$$K \to \left(L^K - \lambda\right)^{-1}, \, Im \lambda > 0$$
являются при каждом фиксированном $\lambda$ гомеоморфизмами (см. замечание 4).

Перейдем к доказательствам сформулированных результатов.
Доказательству Основной леммы предпошлем две леммы, являющиеся частными случаями соответствующих
 утверждений для общих квазидифференциальных выражений (см. [\ref{EM}]).

\textbf{Лемма 6.}
\emph{Пусть} $y, z \in \text{Dom}(L_{\text{max}}).$\emph{ Тогда}
$$
\int\limits_a^b \left(D^{[m]}y\cdot\overline z  -
y\cdot\overline{D^{[m]}z} \right)dx =
\sum\limits_{k = 1}^{m}
(-1)^{k - 1}{D^{[m - k]}y\cdot\overline {D^{[k -1]}z}}\left|_{x = a}^{x = b}
\right.
$$

\textbf{Лемма 7.}
\emph{Для произвольных
наборов комплексных чисел} $ \{\alpha _0 ,\alpha _1, ..., \alpha _{m - 1} \}$, $\{ \beta _0
,\beta _1, ..., \beta _{m - 1} \}$ \emph{существует функция }${y \in \text{Dom}(L_{\text{max}})}$\emph{ такая, что}
$$ D^{[k]}y(a) = \alpha _k , \quad D^{[k]}y(b) = \beta _k, \quad k = 0,1, ..., m - 1.$$

\textbf{\emph{Доказательство Основной леммы.}}
Достаточно показать, что тройка  $(\mathbb{C}^{m}, \Gamma_1, \Gamma_2)$ удовлетворяет условиям $1)$ и $2)$
определения ПГЗ с $\mathcal{H} = L_2$.
Согласно теореме 1, $L^*_{\text{min}} = L_{\text{max}}$.
Согласно лемме 6
\[\left( {L_{\text{max}}y,z} \right) - \left( y,L_{\text{max}}z\right)  =
i^m\sum\limits_{k = 1}^{m}(-1)^{k - 1}{D^{[m - k]}y\cdot\overline {D^{[k -1]}z}}\left|_{x = a}^{x = b}
\right..\]

Однако легко подсчитать, что в случае $m = 2n$
 $$\left( {\Gamma _1 y,\Gamma _2 z} \right) = i^{2n}\sum\limits_{k = 1}^{n}
(-1)^{k - 1}{D^{[2n - k]}y\cdot\overline {D^{[k -1]}z}}\left|_{x = a}^{x = b}\right.,$$
\noindent и
\[\left({\Gamma _2 y,\Gamma _1 z} \right) = i^{2n}\sum\limits_{k = n + 1}^{2n}
(-1)^{k}{D^{[2n - k]}y\cdot\overline {D^{[k -1]}z}}\left|_{x = a}^{x = b}\right..\]

Это показывает, что выполняется условие $1)$. Выполнение условия $2)$ следует непосредственно из леммы 7.

В случае $m = 2n + 1$ введем обозначения :
$$\Gamma _1 =: \left(\Gamma _{1a}, \Gamma _{1b}, \Gamma _{1ab}\right),$$
$$\Gamma _2 =: \left(\Gamma _{2a}, \Gamma _{2b}, \Gamma _{2ab}\right),$$
где
\begin{align*}
&\Gamma _{1a} = i^{2n + 1}\left( - D^{[2n]}y(a), ... , (-1)^{n} D^{[n + 1]}y(a)\right),\\
&\Gamma _{1b} = i^{2n + 1}\left(D^{[2n]}y(b), ... , (-1)^{n + 1} D^{[n + 1]}y(a)\right), \\
&\Gamma _{1ab} = i^{2n + 1}\left(\alpha D^{[n]}y(b) + \beta D^{[n]}y(a)\right),\\
&\Gamma _{2a} = \left(D^{[0]}y(a), ..., D^{[n - 1]}y(a)\right),\\
&\Gamma _{2b} = \left(D^{[0]}y(b), ..., D^{[n - 1]}y(b)\right),\\
&\Gamma _{2ab} = \gamma D^{[n]}y(b) + \delta D^{[n]}y(a).
\end{align*}

Тогда легко проверить, что

$\left( {\Gamma _{1a} y,\Gamma _{2a} z} \right)  = i^{2n + 1}
    \sum\limits_{k = 1}^{n}(-1)^{k - 1}{D^{[2n - k]}y(a)\cdot\overline {D^{[k -1]}z(a)}},$

$\left({\Gamma _{2a} y, \Gamma _{1a} z} \right) = i^{2n + 1}
    \sum\limits_{k = n + 2}^{2n + 1}(-1)^{k}{D^{[2n - k]}y(a)\cdot\overline {D^{[k -1]}z(a)}},$

$\left( {\Gamma _{1b} y,\Gamma _{2b} z} \right)  = i^{2n + 1}
    \sum\limits_{k = 1}^{n}(-1)^{k - 1}{D^{[2n - k]}y(b)\cdot\overline {D^{[k -1]}z(b)}},$

$\left({\Gamma _{2b} y, \Gamma _{1b} z} \right) = i^{2n + 1}
    \sum\limits_{k = n + 2}^{2n + 1}(-1)^{k}{D^{[2n - k]}y(b)\cdot\overline {D^{[k -1]}z(b)}},$

\begin{align*}
\left( {\Gamma _{1ab} y,\Gamma _{2ab} z} \right) -& \left( {\Gamma _{2ab} y,\Gamma _{1ab} z} \right) =\\
= i^{2n + 1}&(-1)^{n} \left(D^{[n]}y(b)\cdot\overline {D^{[n]}z(b)} - D^{[n]}y(a)\cdot\overline {D^{[n]}z(a)}\right).
\end{align*}

Из приведенных соотношений следует, что выполняется условие 1) определения 3.
Выполнение условия $2)$ следует из леммы 7 и последнего из соотношений (\ref{PGZ coef}).

\textbf{\emph{Доказательство теорем 4, 5 и 6.}}
Утверждения теорем следуют из Основной леммы и теоремы 1.6 гл. 3 монографии [\ref{Gorbachuk}]
для ПГЗ абстрактного симметрического оператора.

%%%%%%%%%%%%%%%%%%%%%%%%%%%%%%%%%%%%%%%%%%%%%%%%%%%%%%%%%

\section{Обобщенные резольвенты.}

%\textbf{5. Обобщенные резольвенты.}
Напомним известное

\textbf{Определение 5.}
\emph{Обобщенной резольвентой замкнутого симметрического оператора} $L$\emph{ называют
операторную функцию }$R_\lambda$\emph{ комплексного параметра }$\lambda \in \mathbb{C} \backslash \mathbb{R},$
\emph{допускающую представление вида}
$$ R_\lambda f = P^+ \left( L^+ - \lambda I^+\right)^{- 1}f, \quad f \in \mathcal{H},$$
\emph{где }$L^+$ \emph{-- какое-либо самосопряженное расширение оператора }$L$\emph{ с выходом, вообще говоря,
в более широкое, чем }$\mathcal{H}$, \emph{пространство} $\mathcal{H}^+,$
$I^+$\emph{ -- единичный оператор в }$\mathcal{H}^+,$
$P^+$ \emph{-- оператор ортогонального проектирования }$\mathcal{H}^+$ \emph{на} $\mathcal{H}.$

Операторная функция $R_\lambda$ $(\text{Im} \lambda \neq 0)$ является обобщенной
резольвентой симметрического оператора $L$ тогда и только тогда,
когда
$$\left( R_\lambda f, g \right)_\mathcal{H} = \int_{-\infty}^{+\infty}\frac{d\left(F_\mu f, g\right)}{\mu - \lambda},
\quad f, g \in \mathcal{H}, $$
где $F_\mu$ -- обобщенная спектральная функция оператора $L$.
Это означает, что операторная функция $F_\mu$, $\mu \in \mathbb{R},$
обладает следующими свойствами [\ref{Ahiezer}]:

$1^0.$ При $\mu_2 > \mu_1$ разность $F_{\mu_2} - F_{\mu_1}$ является
ограниченным неотрицательным оператором;

$2^0.$ $F_{\mu +} = F_\mu$ при всех вещественных $\mu$;

$3^0.$ При любом $x \in \mathcal{H}$ $$ \lim\limits_{\mu \rightarrow
- \infty}^{}||F_\mu x ||_\mathcal{H} = 0, \quad \lim\limits_{\mu
\rightarrow + \infty}^{} ||{F_\mu x - x} ||_\mathcal{H} = 0.$$

Следующий результат принадлежит В. М. Бруку [\ref{Brook}]. %(см. Т. 1 и Зам. 1).

Пусть $H$ --- вспомогательное сепарабельное гильбертово пространство.
Символом $\left\{X, X'\right\}$ будем обозначать упорядоченную пару из $X, X' \in H$.
Пары $\left\{X, X'\right\}$ рассматриваются как элементы пространства $H \oplus H$.
Предположим, что существует линейный оператор $\gamma$, отображающий область определения $\operatorname{Dom}(L^*)$
сопряженного к $L$ оператора $L^*$ на $H \oplus H$,
такой, что имеет место равенство
$$
(L^*x, y) - (x, L^*y) = (X', Y)_H - (X, Y')_H,
$$
где $x, y \in \operatorname{Dom}(L)$, $\left\{X, X'\right\} = \gamma x$, $\left\{Y, Y'\right\} = \gamma y$.

\textbf{Теорема 7.}[\ref{Brook}]
\emph{Существует взаимно однозначное соответствие между обобщенными
резольвентами оператора $L$и краевыми задачами}
\[
L^*y = \lambda y + h,
\]
\[
\left( {K(\lambda) - I} \right)Y' \mp i\left( {K(\lambda) + I} \right)Y = 0
\]
\emph{где $\{Y,Y'\} = \gamma y$, $h \in \mathcal{H}$, $\lambda$ --- комплексное число и знак \glqq\,$+$\grqq  \,
в краеовом условии берется для значений
$\lambda$ из верхней полуплоскости,
а \glqq\,$-$\grqq \, для $\lambda$ из нижней полуплоскости.
%$\operatorname{Im}(\lambda) > 0$,
$K(\lambda)$ --- заданная регулярная в верхней полуплоскости операторная функция в $H$, такая, что
$\|K(\lambda)\|_H \leq 1$;
для значений $\lambda$ из нижней полуплоскости $K(\lambda)$ определяется как $K^*(\overline{\lambda})$.}

\emph{Каждое решение задачи определяет обобщенную резольвенту оператора $L$ и обратно, каждая обощенная резольвента оператора
$L$ определяется решением данной задачи.}

Эта теорема дает возможность описать все обобщенные резольвенты симметрического оператора $L_{\operatorname{min}}$
вне вещественной оси.

Параметрическое \emph{внутреннее} описание всех обобщенных резольвент
 оператора $L_{\text{min}}$
дает

\textbf{Теорема 8.}
$1)$ \emph{Каждая обобщенная резольвента $R_\lambda$ оператора $L_{\operatorname{min}}$
в полуплоскости $\operatorname{Im}\lambda < 0$ задается формулой $R_\lambda h = y$,
где $y$ --- решение краевой задачи вида}
%де
%Існує взаємно однозначна відповідність між узагальненими
%резольвентами оператора $L_{\operatorname{min}}$ і крайовими задачами
$$ l(y) = \lambda y + h,$$
$$
 \left( {K(\lambda) - I} \right)\Gamma _{[1]} f + i\left( {K(\lambda) + I} \right)\Gamma _{[2]} f = 0.
$$
\emph{Тут $h(x) \in L_2\left(\mathcal{J},\mathbb{C}\right)$ і
$K(\lambda)$ --- регулярная в нижней полуплоскости операторная функция в пространстве $\mathbb{C}^{2}$ такая,
что $||K(\lambda)|| \leq 1$.}

$2)$ \emph{В полуплоскости $\operatorname{Im}\lambda > 0$
каждая обобщенная резольвента оператора $L_{\operatorname{min}}$ задается формулой $R_\lambda h = y$,
где $y$ --- решение краевой задачи вида}
%де
%Існує взаємно однозначна відповідність між узагальненими
%резольвентами оператора $L_{\operatorname{min}}$ і крайовими задачами
$$ l(y) = \lambda y + h,$$
$$
 \left( {K(\lambda) - I} \right)\Gamma _{[1]} f - i\left( {K(\lambda) + I} \right)\Gamma _{[2]} f = 0.
$$
\emph{Тут $h(x) \in L_2\left(\mathcal{J},\mathbb{C}\right)$ и
$K(\lambda)$ --- регулярная в верхней полуплоскости операторная функция в пространстве $\mathbb{C}^{2}$ такая,
что $||K(\lambda)|| \leq 1$.}

\emph{Эта параметризация обобщенных резольвент операторными функциями $K$ является биективной.}

\textbf{\emph{Доказательство.}}
В силу Основной леммы вспомогательное сепарабельное гильбертово пространство $\mathbb{C}^m$ и оператор
$\gamma y = \{\Gamma_{[1]}y, \Gamma_{[2]}y\}$, отображающий $\operatorname{Dom}(L_{\operatorname{min}})$ на $\mathbb{C}^m\oplus\mathbb{C}^m$, удовлетворяют условиям теоремы 7.

Таким образом, утверждение теоремы 8 вытекает из теоремы 7.

%%%%%%%%%%%%%   СПИСОК ЛИТЕРАТУРЫ %%%%%%%%%%%%%%%
\vskip 3.5mm

\footnotesize

\renewcommand{\labelenumi}{[\theenumi]}
\begin{enumerate}
\item\label{Z}
{\it Zettl~A.}
Sturm--Liouville Theory.~--- Providence:
{American Mathematical Society}, {2005}.~--- 328 p.

\item\label{S-Sh}
{\it Савчук~А.М., Шкаликов~А.А.}
Операторы Штурма--Лиувилля с сингулярными потенциалами~//
Мат. заметки.~--- 1999.~--- \textbf{66}, №~6.~--- С. 897--912.

\item\label{Albeverio}
{\it Albeverio~S., Gestezy~F., Hoegh-Krohn~R., Holden~H.}
Solvable models in quantum mechanics.~--- {New York}: {Springer-Verlag},	
{1988}.~--- 452 p.	

\item\label{GM St-L}
 {\it Goriunov~A.S., Mikhailets~V.A.}  %authors
Regularization of singular Sturm--Liouville equations~//
Methods Funct. Anal. Topology.~--- 2010.~--- №~2.~--- С. 120--130.

\item\label{Shin}
{\it Шин Д.}   %authors
О квазидифференциальных операторах в гильбертовом пространстве~//
Мат. сборник.~--- 1943.~--- \textbf{13(55)}, №1.~--- С. 39--70.

\item\label{Zettl}
{\it Zettl~A.}
Formally self-adjoint quasi-differential operators~//
{Rocky Mountain J. Math.}~--- 1975.~--- \textbf{5}, №~3.~--- P. 453--474.

\item\label{EM}
{\it Everitt~W.N., Markus~L.}
Boundary Value Problems and Symplectic Algebra for Ordinary
Differential and Quasi-differential Operators.~--- {Providence}: {American Mathematical Society},
 {1999}.~--- 187 p.

\item\label{K}
{\it Като~Т.}
Теория возмущений линейных операторов.~--- М.: Мир, 1972.~--- 740 с.

\item\label{MR1}
 {\it Михайлец~В.А., Рева~Н.В.}
Обобщения теоремы Кигурадзе о корректности линейных краевых задач~//
Доп. НАН України.~--- 2008.~--- №~9.~--- С. 23--27.

\item\label{MR2}
 {\it Михайлец~В.А., Рева~Н.В.}
Непрерывность по параметру решений общих краевых задач~//
Зб. праць Ін-ту математики.~--- 2008.~--- \textbf{5}, №~1.~--- С. 227--239.

\item\label{Levin}
 {\it Левин~А.Ю.}
Предельный переход для несингулярных систем $\dot{X} = A_n(t)X$~//
Докл. АН СССР.~--- 1967.~--- \textbf{176}, №~4.~--- С. 774--777.

\item\label{Gorbachuk}
{Горбачук~В.И., Горбачук~М.Л.}
Граничные задачи для дифференциально-операторных уравнений.~--- К.: Наук. думка, 1984.~--- 284 с.

 \item\label{Phil}
 {Phillips~R.S.}
Dissipative operators and hyperbolic systems of partial differential equations~//
{Trans. Amer. Math. Soc.}~--- 1959.~--- \textbf{90}.~--- P.~193--254.

\item\label{Ahiezer}
{Ахиезер~Н.И., Глазман~И.М.}
Теория линейных операторов в гильбертовом пространстве.~--- М.: Наука, 1966.~--- 544 с.

\item\label{Brook}
{Брук~В.М.}
Об одном классе краевых задач со спектральным параметром в граничном условии~//
Мат. сборник.~--- 1976.~---\textbf{100(142)}, {№~2(6)}.~--- С. 210--216.
\end{enumerate}
\end{document}